\documentclass[reqno, letterpaper, 11pt]{amsart}

\usepackage[T1]{fontenc}
\usepackage[utf8]{inputenc}
\usepackage{lineno}
\usepackage{amssymb}
\usepackage{graphicx}
\usepackage{latexsym}
\usepackage{xcolor}
\usepackage{microtype}
\usepackage{amsmath,amsthm,mathtools}
\usepackage{aliascnt}
\usepackage{enumitem}
\usepackage[bookmarks]{hyperref}
\hypersetup{
  unicode,
  colorlinks,
  linkcolor={red!60!black},
  citecolor={green!60!black},
  urlcolor={blue!60!black},
  pdftitle={Ray inflations of omega1-trees: subdivisions, subgraph rigidity, and ends of degree aleph1},
  pdfauthor={Leandro Aurichi, Gabriel Fernandes, and Paulo Magalhaes Junior},
  pdfsubject={Ray inflations of sparse graphs on omega1-trees},
  pdfkeywords={omega1-tree, Aronszajn tree, almost-Suslin tree, ray inflation, graph end, subdivision}
}

\usepackage[nameinlink, capitalise, noabbrev]{cleveref}
\usepackage{geometry}
\usepackage[english]{babel}
\usepackage{mathptmx}
\usepackage[font=footnotesize]{caption}
\numberwithin{equation}{section}

\setenumerate{label={\normalfont (\roman*)},itemsep=0pt}
\setlist[itemize]{leftmargin=2em,itemsep=2pt,topsep=4pt}

\newtheorem{theorem}{Theorem}[section]

\newaliascnt{lemma}{theorem}
\newtheorem{lemma}[lemma]{Lemma}
\aliascntresetthe{lemma}

\newaliascnt{proposition}{theorem}
\newtheorem{proposition}[proposition]{Proposition}
\aliascntresetthe{proposition}

\newaliascnt{corollary}{theorem}
\newtheorem{corollary}[corollary]{Corollary}
\aliascntresetthe{corollary}

\newaliascnt{conjecture}{theorem}
\newtheorem{conjecture}[conjecture]{Conjecture}
\aliascntresetthe{conjecture}

\newaliascnt{question}{theorem}
\newtheorem{question}[question]{Question}
\aliascntresetthe{question}

\theoremstyle{definition}
\newaliascnt{definition}{theorem}
\newtheorem{definition}[definition]{Definition}
\aliascntresetthe{definition}

\theoremstyle{remark}
\newaliascnt{remark}{theorem}
\newtheorem{remark}[remark]{Remark}
\aliascntresetthe{remark}

\theoremstyle{plain}
\newaliascnt{theorem1}{theorem}
\newtheorem*{theorem1}{Theorem \ref{thm:intro-nonsubdivision}}
\aliascntresetthe{theorem1}

\theoremstyle{plain}
\newaliascnt{theorem2}{theorem}
\newtheorem*{theorem2}{Theorem \ref{thm:diamond-example}}
\aliascntresetthe{theorem2}

\theoremstyle{plain}
\newaliascnt{theorem3}{theorem}
\newtheorem*{theorem3}{Theorem \ref{thm:intro-rigidity}}
\aliascntresetthe{theorem3}

\theoremstyle{plain}
\newaliascnt{corollary4}{theorem}
\newtheorem*{corollary4}{Corollary \ref{cor:intro-almost-suslin}}
\aliascntresetthe{corollary4}

\crefname{theorem}{Theorem}{Theorems}
\Crefname{theorem}{Theorem}{Theorems}
\crefname{lemma}{Lemma}{Lemmas}
\Crefname{lemma}{Lemma}{Lemmas}
\crefname{proposition}{Proposition}{Propositions}
\Crefname{proposition}{Proposition}{Propositions}
\crefname{corollary}{Corollary}{Corollaries}
\Crefname{corollary}{Corollary}{Corollaries}
\crefname{definition}{Definition}{Definitions}
\Crefname{definition}{Definition}{Definitions}
\crefname{conjecture}{Conjecture}{Conjectures}
\Crefname{conjecture}{Conjecture}{Conjectures}
\crefname{question}{Question}{Questions}
\Crefname{question}{Question}{Questions}
\crefname{remark}{Remark}{Remarks}
\Crefname{remark}{Remark}{Remarks}

\newcommand{\treeht}{\operatorname{ht}}
\newcommand{\treepred}{\operatorname{pred}}
\newcommand{\Succ}{\operatorname{Succ}}
\newcommand{\acc}{\operatorname{acc}}
\newcommand{\degop}{\operatorname{deg}}
\newcommand{\downset}[1]{\lceil #1\rceil}
\newcommand{\upset}[1]{\lfloor #1\rfloor}
\newcommand{\sdownset}[1]{\mathring{\lceil #1\rceil}}
\newcommand{\supsetopen}[1]{\mathring{\lfloor #1\rfloor}}
\newcommand{\restr}{\mathbin{\upharpoonright}}
\newcommand{\rayinfl}[1]{#1\mathbin{\#}\mathbb{N}}
\newcommand{\HC}{\mathrm{HC}}
\newcommand{\HCstar}{\mathrm{HC}^{*}}

\title{\texorpdfstring{Ray inflations of $\omega_1$-trees and ends of degree $\aleph_1$}{Ray inflations of omega1-trees and ends of degree aleph1}}

\author{Leandro Aurichi}
\address{Instituto de Ci\^encias Matem\'aticas e de Computa\c{c}\~ao, Universidade de S\~ao Paulo, Avenida Trabalhador S\~ao-Carlense 400, 13566-590 S\~ao Carlos, SP, Brazil}
\email{aurichi@icmc.usp.br}

\author{Gabriel Fernandes}
\address{Instituto de Ci\^encias Matem\'aticas e de Computa\c{c}\~ao, Universidade de S\~ao Paulo, Avenida Trabalhador S\~ao-Carlense 400, 13566-590 S\~ao Carlos, SP, Brazil}
\email{fernandes@icmc.usp.br}

\author{Paulo Magalh\~aes J\'unior}
\address{Instituto de Ci\^encias Matem\'aticas e de Computa\c{c}\~ao, Universidade de S\~ao Paulo, Avenida Trabalhador S\~ao-Carlense 400, 13566-590 S\~ao Carlos, SP, Brazil}
\email{pjr.mat@usp.br}


\begin{document}

\begin{abstract}
We prove the followings result for ray inflations of sparse graphs on $\omega_1$-trees. First, let $T$ and $S$ be pruned $\omega_1$-trees, let $G_T$ be a sparse $T$-graph, and let $G_S$ be a sparse $S$-graph with uniformly finite adhesion. If $T$ is not special, then no subdivision of $\rayinfl{G_S}$ is isomorphic to $\rayinfl{G_T}$. Second, if $T$ is almost-Suslin and $S$ is special, then $\rayinfl{G_T}$ contains no subgraph isomorphic to $\rayinfl{G_S}$, for arbitrary choices of the sparse graphs. Consequently, under $\diamondsuit_{\omega_1}$, this gives a counterexample to Halin's end degree conjecture that is not isomorphic to a subdivision of any ray inflation with uniformly finite adhesion. Under $\diamondsuit_{\omega_1}^{*}$, the underlying tree may in addition be chosen almost-Suslin, and the resulting graph contains no subgraph isomorphic to a ray inflation over a special $\omega_1$-tree. 
\end{abstract}

\maketitle

\section{Introduction}\label{sec:introduction}

Originally defined by Halin in \cite{Halin1964}, a ray is a
one-way infinite path; its infinite connected subgraphs are called its \textbf{tails}. We say that two rays in a graph $G$ are \textbf{equivalent} if no finite set of vertices separates them; the
corresponding equivalence classes of rays are the \textbf{ends} of $G$, also defined in \cite{Halin1965}. The set of ends of a graph $G$ is
denoted by $\Omega = \Omega(G)$.

Applications of ends have been found in various areas of mathematics, as surveyed in \cite{MR2800969}. Many aspects of ends have been studied, including their topological, combinatorial, and set-theoretical properties. A key result in the theory of infinite graphs is Halin's grid theorem, which addresses configurations of ends with infinite degree. The degree $\deg(\varepsilon)$ of an end $\varepsilon$ is the supremum of the sizes of collections of pairwise disjoint rays in $\varepsilon$. The supremum is always attained (see, e.g., \cite{Halin1965}). A \textbf{hexagonal grid} is an infinite planar cubic graph that contains several rays in a well-organized pattern. Given a graph $G$, we say that a graph $H$ is a \textbf{subdivision} of $G$ if $H$ is obtained from $G$ by `subdividing' some or all of its edges by drawing new vertices on those edges.

\begin{theorem}[Halin's grid theorem - \cite{Halin1965}]\label{halin1965}
Every graph with an end of infinite degree contains a subdivision of the hexagonal quarter grid whose rays belong to that end.
\end{theorem}

 To extend Theorem \ref{halin1965} to larger cardinalities, Halin introduced the concept of a `ray graph' as follows:

Given a graph $G$, a family of paths in $G$ is called \textbf{independent} if they intersect only at their endpoints. Let $R$ be a set of pairwise disjoint equivalent rays in $G$, and let $P$ be a family of independent paths in $G$ whose endpoints lie precisely on the rays of $R$. The \textbf{ray graph} $H$ associated with $R$ and $P$ is defined as the graph  whose vertex set $V(H)$ consists of the rays in $R$,  two rays $r$,$s$ in $R$ are adjacent if and only if there are infinitely many pairwise disjoint paths in $P$ connecting the rays $r$ and $s$. A \textbf{ray graph for an end} $\varepsilon$ is a connected ray graph in $G$ constructed on a subset of $\varepsilon$ that witnesses the degree of the end.

In \cite{Halin2000}, Halin conjectured the following:

\begin{conjecture}[Halin's Degree Conjecture - \cite{Halin2000}, Conjecture 6.1]
Every graph admits a ray graph for each of its ends.
\end{conjecture}

In \cite{GeschkeEtAl}, it was shown that Halin's Degree Conjecture fails for ends of degree $\aleph_1$, although it holds in other cases, as also proven in \cite{GeschkeEtAl}. For recent developments on Halin's Degree Conjecture for degrees $> \aleph_1$, see \cite{fernandes2025remarkshalinsenddegreeconjecture}.

The counterexample to Halin's Degree Conjecture for $\aleph_1$, presented in \cite{GeschkeEtAl}, involves the ray inflation of a sparse $T$-graph, where $T$ is a special Aronszajn tree. The following question is posed in \cite{GeschkeEtAl}:

\begin{question}[\cite{GeschkeEtAl}, Question 2] \label{CentralQuestion}
Is it true that any end of degree $\aleph_1$ contains either an $\aleph_1$-star of rays or a subdivision of a ray inflation of a sparse $T$-graph, where $T$ is a special Aronszajn tree?
\end{question}

The present paper studies two related but logically distinct comparison problems for ray inflations. The first problem asks whether one ray inflation can be isomorphic to a subdivision of another. We prove the following theorem.

\begin{theorem1}
Let $T$ and $S$ be pruned $\omega_1$-trees. Let $G_T$ be a sparse $T$-graph and let $G_S$ be a sparse $S$-graph with uniformly finite adhesion. If $T$ is not special, then $\rayinfl{G_T}$ is not isomorphic to a subdivision of $\rayinfl{G_S}$.
\end{theorem1}

The characterisation of Kurkofka, Pitz \cite[Theorem~4.6]{KurkofkaPitz} states that an order tree admits a $T$-graph of uniformly finite adhesionif and only if it is special, and that it admits a $T$-graph of finite adhesionif and only if it is semi-special. Their proof gives canonical sparse graphs from the corresponding antichain partitions. Consequently, \cref{thm:intro-nonsubdivision} has the following form.

\begin{theorem2}
Assume $\diamondsuit_{\omega_1}$. There is a graph $I$ with exactly one end $\varepsilon$ such that
\begin{enumerate}
  \item $\degop(\varepsilon)=\aleph_1$ and $\varepsilon$ has no ray graph;
  \item if $S$ is a pruned $\omega_1$-tree and $G_S$ is a sparse $S$-graph with uniformly finite adhesion, then $I$ is not isomorphic to a subdivision of $\rayinfl{G_S}$.
\end{enumerate}
\end{theorem2}

The stronger principle $\diamondsuit_{\omega_1}^{*}$ yields a pruned semi-special non-special almost-Suslin Aronszajn tree. For this choice of the underlying tree, \cref{cor:intro-almost-suslin} also excludes every subgraph isomorphic to a ray inflation over a special $\omega_1$-tree. The prediction principles are used only to establish the existence of the required trees

The second problem asks whether one ray inflation can occur as a subgraph of another without subdividing its edges. The tree-theoretic hypothesis used for this problem is almost-Suslinity. Recall that an $\omega_1$-tree is almost-Suslin if the set of heights of every antichain is nonstationary. We prove the following theorem for arbitrary graph embeddings between ray inflations.

\begin{theorem3}
Let $S$ and $T$ be $\omega_1$-trees, and let $G_S$ and $G_T$ be sparse graphs on $S$ and $T$, respectively. If
\[
  f:\rayinfl{G_S}\longrightarrow\rayinfl{G_T}
\]
is a graph embedding, then there are a club $C\subseteq\omega_1$, consisting of non-zero limit ordinals, and a level-preserving order embedding
\[
  \pi:S\restr C\longrightarrow T\restr C
\]
such that, for every $\alpha\in C$ and $s\in S^\alpha$, the image $f[R_s]$ is a tail of $R_{\pi(s)}$.
\end{theorem3}

\begin{corollary4}
Let $T$ be an almost-Suslin $\omega_1$-tree and let $S$ be a special $\omega_1$-tree. For arbitrary sparse graphs $G_T$ on $T$ and $G_S$ on $S$, the graph $\rayinfl{G_T}$ contains no subgraph isomorphic to $\rayinfl{G_S}$.
\end{corollary4}

The embedding in \cref{cor:intro-almost-suslin} is not required to preserve non-adjacency. Hence the corollary excludes every copy of $\rayinfl{G_S}$ as a subgraph of $\rayinfl{G_T}$. The analogous assertion for a subdivision of $\rayinfl{G_S}$ occurring as a subgraph remains open. \cref{thm:intro-nonsubdivision} excludes an isomorphism from $\rayinfl{G_T}$ onto a subdivision of another ray inflation, whereas \cref{cor:intro-almost-suslin} excludes copies without edge subdivision but does not address topological minors.

Both proofs the Theorems \ref{thm:intro-nonsubdivision} and \ref{thm:intro-rigidity} begin by choosing a club $C\subseteq\omega_1$ on which the rank-initial segments are preserved by the map under consideration. In the proof of \cref{thm:intro-nonsubdivision}, non-specialness yields $\alpha\in C$ and $v\in T^\alpha$ for which infinitely many lower endpoints below $v$ are adjacent to nodes strictly above $v$. Uniformly finite adhesion in the comparison tree bounds the possible images of these lower endpoints by a finite set. In the proof of \cref{thm:intro-rigidity}, write $I_T=\rayinfl{G_T}$. For $s\in S^\alpha$ and every $n\in\omega$, the edge from $(s,n)$ to $(C_n(S,s),n)$ places $f(s,n)$ in the set of vertices of one component of $I_T-(I_T)^{<\alpha}$ that have a neighbour in $(I_T)^{<\alpha}$. The graph induced by this set has $R_{\pi(s)}$ as its unique connected component containing a ray.

The paper is organised as follows. In \cref{sec:preliminaries} we fix the terminology for order trees, the principles $\diamondsuit_{\omega_1}$ and $\diamondsuit_{\omega_1}^{*}$, sparse $T$-graphs, adhesion, ray inflations, ends, and subdivisions. We also record the existence result under $\diamondsuit_{\omega_1}^{*}$ that supplies an almost-Suslin tree with the required additional properties. Section~\ref{sec:club} contains the club arguments used in both proofs of Theorems \ref{thm:intro-nonsubdivision} and \ref{thm:intro-rigidity}. In \cref{sec:whole-graph} we prove \cref{thm:intro-nonsubdivision}. In \cref{sec:subgraph} we prove \cref{thm:intro-rigidity,cor:intro-almost-suslin}. Section~\ref{sec:ends} proves \cref{thm:diamond-example} and introduces $\HCstar(\kappa)$. The open problems are stated in \cref{sec:open}.

\section{Trees, sparse graphs, ray inflations, and ends}\label{sec:preliminaries}

Our set-theoretic terminology follows~\cite{Jech}, and our graph-theoretic terminology follows~\cite{Diestel}. All graphs are simple and undirected; they need not be finite or locally finite. We identify every cardinal with its initial ordinal. The symbol $\omega$ denotes the first infinite ordinal and $\omega_1$ the first uncountable ordinal.

For $A\subseteq\omega_1$, put
\[
  \acc(A):=\{\alpha<\omega_1:\alpha>0\text{ and }\sup(A\cap\alpha)=\alpha\}.
\]
A set $C\subseteq\omega_1$ is a \emph{club} if it is unbounded and $\acc(C)\subseteq C$. A set $A\subseteq\omega_1$ is \emph{stationary} if it meets every club subset of $\omega_1$. The principle $\diamondsuit_{\omega_1}$ asserts that there is a sequence $\langle A_\alpha:\alpha<\omega_1\rangle$ such that $A_\alpha\subseteq\alpha$ for every $\alpha<\omega_1$ and, for every $A\subseteq\omega_1$, the set
\[
  \{\alpha<\omega_1:A\cap\alpha=A_\alpha\}
\]
is stationary.

The principle $\diamondsuit_{\omega_1}^{*}$ asserts that there is a sequence
\[
  \langle\mathcal A_\alpha:\alpha<\omega_1\rangle
\]
such that $\mathcal A_\alpha$ is a countable family of subsets of $\alpha$ for every $\alpha<\omega_1$ and, for every $A\subseteq\omega_1$, there is a club $C_A\subseteq\omega_1$ satisfying
\[
  A\cap\alpha\in\mathcal A_\alpha
  \qquad(\alpha\in C_A).
\]
The formulation above is the form of $\diamondsuit_{\omega_1}^{*}$ used in the existence result below.

\subsection{Order trees}

\begin{definition}\label{def:order-tree}
A partially ordered set $(T,\le_T)$ is an \emph{order tree} if it has a unique minimal element, called the \emph{root}, and
\[
  \downset{t}=\downset{t}_T:=\{t'\in T:t'\le_T t\}
\]
is well ordered by $\le_T$ for every $t\in T$. We write
\[
  \upset{t}=\upset{t}_T:=\{t'\in T:t\le_T t'\},
\]
and abbreviate
\[
  \sdownset{t}:=\downset{t}\setminus\{t\},
  \qquad
  \supsetopen{t}:=\upset{t}\setminus\{t\}.
\]
For $X\subseteq T$, put
\[
  \downset{X}:=\bigcup_{t\in X}\downset{t},
  \qquad
  \upset{X}:=\bigcup_{t\in X}\upset{t}.
\]
The height of $t$ is
\[
  \treeht_T(t):=\operatorname{otp}(\sdownset{t},<_T),
\]
and the height of $T$ is the least ordinal $\alpha$ for which no node has height $\alpha$. The $\alpha$th level is
\[
  T^\alpha:=\{t\in T:\treeht_T(t)=\alpha\}.
\]
For an ordinal $\alpha$, write
\[
  T^{<\alpha}:=\bigcup_{\beta<\alpha}T^\beta,
  \qquad
  T^{\le\alpha}:=\bigcup_{\beta\le\alpha}T^\beta.
\]
If $\kappa$ is a regular uncountable cardinal, a \emph{$\kappa$-tree} is an order tree of height $\kappa$ whose levels have cardinality less than $\kappa$.
\end{definition}

A maximal chain is a \emph{branch}. A node is \emph{terminal} if it has no strict extension, and a tree is \emph{pruned} if it has no terminal nodes. If $t<_Tt'$ and there is no node strictly between them, then $t'$ is a \emph{successor node}, $t$ is its predecessor, and we write $t=\treepred_T(t')$. A non-root node with no predecessor is a \emph{limit node}. We denote the set of successor nodes by $\Succ(T)$. If $\beta\le\treeht_T(t)$, the unique node of $T^\beta$ below $t$ is denoted by $t\restr\beta$.

Two nodes are \emph{compatible} if they are comparable; otherwise they are \emph{incompatible}. An \emph{antichain} is a set of pairwise incompatible nodes.

\begin{definition}\label{def:tree-properties}
An order tree is
\begin{enumerate}
  \item \emph{special} if it is the union of countably many antichains;
  \item \emph{semi-special} if its successor nodes are the union of countably many antichains.
\end{enumerate}
An \emph{Aronszajn tree} is an $\omega_1$-tree with no uncountable branch. An $\omega_1$-tree is \emph{almost-Suslin} if, for every antichain $A\subseteq T$, the set
\[
  \{\treeht_T(a):a\in A\}
\]
is nonstationary in $\omega_1$.
\end{definition}

The almost-Suslin terminology goes back to Devlin and Shelah~\cite{DevlinShelah}. Almost-Suslinity alone does not imply that a tree is Aronszajn: the chain $\omega_1$ is almost-Suslin. Every semi-special $\omega_1$-tree is Aronszajn, because an uncountable branch contains uncountably many successor nodes, while each antichain meets a branch in at most one node.

An order tree $T$ is \emph{$\mathbb R$-embeddable} if there is a map $e:T\to\mathbb R$ such that
\[
  s<_Tt\quad\Longrightarrow\quad e(s)<e(t).
\]

\begin{theorem}[Devlin and Shelah - \cite{DevlinShelah}]\label{thm:diamond-star-tree}
Assume $\diamondsuit_{\omega_1}^{*}$. There is a pruned $\mathbb R$-embeddable almost-Suslin Aronszajn $\omega_1$-tree.
\end{theorem}

\begin{corollary}\label{cor:diamond-star-tree}
Assume $\diamondsuit_{\omega_1}^{*}$. There is a pruned semi-special non-special almost-Suslin Aronszajn $\omega_1$-tree.
\end{corollary}

\begin{proof}
Let $T$ be given by \cref{thm:diamond-star-tree}, and fix a strict order embedding $e:T\to\mathbb R$. For every successor node $t$, choose a rational number $q(t)$ such that
\[
  e(\treepred_T(t))<q(t)<e(t).
\]
For each $q\in\mathbb Q$, the set $\{t\in\Succ(T):q(t)=q\}$ is an antichain. Indeed, if $s<_Tt$ are successor nodes with $q(s)=q(t)$, then $s\le_T\treepred_T(t)$, and hence
\[
  q(s)<e(s)\le e(\treepred_T(t))<q(t),
\]
a contradiction. Thus $T$ is semi-special.

It remains to show that $T$ is not special. Suppose that $T=\bigcup_{n<\omega}A_n$, with each $A_n$ an antichain. By almost-Suslinity, each set
\[
  \{\treeht_T(a):a\in A_n\}
\]
is nonstationary. Their countable union is therefore nonstationary. On the other hand, $T^\alpha\ne\varnothing$ for every $\alpha<\omega_1$, so the union contains every ordinal below $\omega_1$ and is equal to $\omega_1$. This is impossible because $\omega_1$ is stationary.
\end{proof}

\begin{remark}\label{rem:diamond-star-role}
The principles $\diamondsuit_{\omega_1}$ and $\diamondsuit_{\omega_1}^{*}$ are not used in the graph-theoretic arguments. The principle $\diamondsuit_{\omega_1}$ is used to obtain a semi-special non-special tree for the example in \cref{thm:diamond-example}. The principle $\diamondsuit_{\omega_1}^{*}$ is used when the underlying tree is additionally required to be almost-Suslin, as in \cref{thm:diamond-star-example}.
\end{remark}

\subsection{\texorpdfstring{Sparse $T$-graphs and adhesion}{Sparse T-graphs and adhesion}}

For a graph $G$ with vertex set an order tree $T$, put
\[
  N_G^-(t):=N_G(t)\cap\sdownset{t}.
\]
For $X\subseteq V(G)$, let
\[
  N_G(X):=\{v\in V(G)\setminus X:\text{some }x\in X\text{ is adjacent to }v\}
\]
be its external neighbourhood.

\begin{definition}\label{def:T-graph}
An order tree $T$ is \emph{normal in a graph $G$} if $V(G)=T$ and the endvertices of every edge of $G$ are compatible in $T$. A graph $G$ is a \emph{$T$-graph} if $T$ is normal in $G$ and $N_G^-(t)$ is cofinal in $\sdownset{t}$ for every non-root node $t$.
\end{definition}

\begin{definition}\label{def:sparse}
Let $T$ be an order tree of height at most $\omega_1$. A \emph{sparse $T$-graph} is specified by choosing, for every limit node $v$, a strictly increasing sequence
\[
  C(T,v)=\langle C_n(T,v):n\in\omega\rangle
\]
cofinal below $v$, and then setting $V(G)=T$ and
\[
\begin{split}
 E(G)={}&\bigl\{\{\treepred_T(v),v\}:v\in\Succ(T)\bigr\}\\
 &{}\cup\bigl\{\{C_n(T,v),v\}:v\text{ is a limit node},\ n\in\omega\bigr\}.
\end{split}
\]
The family of sequences $\langle C(T,v):v\text{ is a limit node}\rangle$ is part of the data of $G$.
\end{definition}

Thus an order tree does not determine a unique sparse graph unless the cofinal lower-neighbour sequences are fixed.

\begin{definition}\label{def:adhesion}
Let $G$ be a $T$-graph.
\begin{enumerate}
  \item The graph $G$ has \emph{finite adhesion} if, for every limit ordinal $\sigma$, each component of $G-T^{\le\sigma}$ has finite external neighbourhood.
  \item The graph $G$ has \emph{uniformly finite adhesion} if, for every limit node $t\in T$, there is a finite set $S_t\subseteq\sdownset{t}$ such that
  \[
    N_G^-(t')\cap\sdownset{t}\subseteq S_t
    \qquad\text{for every }t'>_Tt.
  \]
\end{enumerate}
\end{definition}

We use the following characterisation.

\begin{theorem}[Kurkofka and Pitz~\protect-\cite{KurkofkaPitz}]\label{thm:KP}
Let $T$ be an order tree.
\begin{enumerate}
  \item There is a $T$-graph with uniformly finite adhesion if and only if $T$ is special.
  \item There is a $T$-graph with finite adhesion if and only if $T$ is semi-special.
\end{enumerate}
\end{theorem}

We shall also use the equivalent formulation that a $T$-graph has finite adhesion if and only if $N_G(\upset{t})$ is finite for every successor node $t$; see~\cite[Lemma~4.2]{KurkofkaPitz}.

The next standard construction fixes the sparse graphs associated with antichain decompositions.

\begin{proposition}[Diestel and Leader~\protect-\cite{DiestelLeader}]\label{prop:standard-graph}
Let $T$ be an $\omega_1$-tree.
\begin{enumerate}
  \item A partition $\mathcal A=\langle A_j:j\in\omega\rangle$ of $\Succ(T)$ into antichains determines a sparse $T$-graph $G_{\mathcal A}$ with finite adhesion.
  \item A partition $\mathcal A=\langle A_j:j\in\omega\rangle$ of $T$ into antichains determines a sparse $T$-graph $G_{\mathcal A}$ with uniformly finite adhesion.
\end{enumerate}
\end{proposition}

\begin{proof}
For every limit node $v$, recursively define a sequence $C(T,v)=\langle C_n(T,v):n\in\omega\rangle$. Let
\[
  \ell_0(v):=\min\{j:A_j\cap\sdownset{v}\ne\varnothing\},
\]
and let $C_0(T,v)$ be the unique element of $A_{\ell_0(v)}\cap\sdownset{v}$. Having chosen $C_n(T,v)$, let
\[
  \ell_{n+1}(v):=\min\{j:A_j\cap\{x:C_n(T,v)<_Tx<_Tv\}\ne\varnothing\},
\]
and let $C_{n+1}(T,v)$ be the unique element of $A_{\ell_{n+1}(v)}\cap\{x:C_n(T,v)<_Tx<_Tv\}$. The intersections are singletons because each $A_j$ is an antichain and $\sdownset{v}$ is a chain. They are non-empty because the successor nodes below a limit node are cofinal below it.

The resulting sequence is cofinal below $v$. Suppose otherwise. Choose a successor node $x<_Tv$ above every $C_n(T,v)$, and let $j$ be the unique index such that $x\in A_j$. At each stage, the node $x$ belongs to the interval from which the next term is selected, so the selected antichain index is at most $j$. Since each of the antichains $A_0,\ldots,A_j$ meets the chain $\sdownset{v}$ in at most one node, at most $j+1$ terms can be selected. This contradicts that the recursion defines $C_n(T,v)$ for every $n\in\omega$.

For (i), let $x\in A_j$ be a successor node. An edge from $\upset{x}$ to its complement is either the edge from $x$ to $\treepred_T(x)$ or has the form $uC_n(T,u)$, where $u>_Tx$ is a limit node and $C_n(T,u)<_Tx$. In the latter case, $C_n(T,u)<_Tx<_Tu$. Thus, at the stage at which $C_n(T,u)$ was selected, the antichain $A_j$ contained an eligible node, namely $x$. By the choice of the least eligible index, $C_n(T,u)\in A_i$ for some $i\le j$. Hence
\[
  N_{G_{\mathcal A}}(\upset{x})
  \subseteq
  \{\treepred_T(x)\}\cup
  \left(\sdownset{x}\cap\bigcup_{i\le j}A_i\right),
\]
which is finite.

For (ii), let $t\in A_j$ be a limit node. If $u>_Tt$ and $C_n(T,u)<_Tt$, then $C_n(T,u)<_Tt<_Tu$. Hence $t\in A_j$ was an eligible node when $C_n(T,u)$ was selected, and the least eligible antichain index was at most $j$. Therefore $C_n(T,u)\in A_i$ for some $i\le j$. Thus
\[
  S_t:=\sdownset{t}\cap\bigcup_{i\le j}A_i
\]
is a finite witness to uniformly finite adhesion at $t$.
\end{proof}

When $\mathcal A$ partitions $\Succ(T)$, we call $G_{\mathcal A}$ the \emph{standard finite-adhesion sparse graph} associated with $(T,\mathcal A)$. When $\mathcal A$ partitions all of $T$, we call it the \emph{standard uniformly finite-adhesion sparse graph} associated with $(T,\mathcal A)$.

\subsection{Ray inflations and subdivisions}

\begin{definition}\label{def:inflation}
Let $G$ be a sparse $T$-graph, with its lower-neighbour sequences fixed. The \emph{ray inflation} $\rayinfl{G}$ has vertex set $T\times\omega$ and the following edges:
\begin{enumerate}
  \item $\{(t,n),(t,n+1)\}$ for every $t\in T$ and $n\in\omega$;
  \item $\{(t,n),(\treepred_T(t),n)\}$ for every successor node $t$ and $n\in\omega$;
  \item $\{(v,n),(C_n(T,v),n)\}$ for every limit node $v$ and $n\in\omega$.
\end{enumerate}
For $t\in T$, the ray induced by $\{t\}\times\omega$ is denoted by $R_t$ and called the \emph{horizontal ray} at $t$.
\end{definition}

We write $\rayinfl{G}$ rather than $T\mathbin{\#}\mathbb{N}$ because the inflation depends on the chosen sparse $T$-graph, not only on the order tree.

Let $I=\rayinfl{G}$, where $G$ is a sparse $T$-graph. Define
\[
  \rho_T(t,n):=\treeht_T(t),
  \qquad
  I^{<\alpha}:=\{x\in V(I):\rho_T(x)<\alpha\}.
\]
We shall repeatedly use the following description of the components obtained after deleting a rank-initial segment.

\begin{lemma}[Geschke, Kurkofka, Melcher, and Pitz~\protect-\cite{GeschkeEtAl}]\label{lem:components-inflation}
For every $\alpha<\omega_1$, the map
\[
  T^\alpha\longrightarrow
  \{\text{components of }I-I^{<\alpha}\},
  \qquad
  t\longmapsto I[\upset{t}\times\omega],
\]
is a bijection.
\end{lemma}

For $t\in T^\alpha$, denote this component by
\[
  K_I(t,\alpha):=I[\upset{t}\times\omega].
\]

A \emph{subdivision} of a graph $H$ is obtained by replacing each edge $xy$ by a finite $x$--$y$ path so that the replacement paths are internally disjoint. Vertices corresponding to $V(H)$ are called \emph{original vertices}; the other vertices are \emph{subdivision vertices}. Every subdivision vertex has degree two.

Let $J$ be a subdivision of $I$. If $z$ is the original vertex corresponding to $(t,n)$, put $\rho_J(z)=\treeht_T(t)$. If $z$ is an internal vertex of the path replacing an edge $xy$ of $I$, put
\[
  \rho_J(z):=\max\{\rho_T(x),\rho_T(y)\}.
\]
Write $J^{<\alpha}:=\{z:\rho_J(z)<\alpha\}$.

\begin{lemma}\label{lem:components-subdivision}
For every $\alpha<\omega_1$, the components of $J-J^{<\alpha}$ are naturally indexed by $T^\alpha$. The component indexed by $t\in T^\alpha$, denoted by $K_J(t,\alpha)$, contains every original vertex $(u,n)$ with $u\ge_Tt$.
\end{lemma}

\begin{proof}
If both endpoints of an edge of $I$ have rank at least $\alpha$, then every vertex of its replacement path remains in $J-J^{<\alpha}$, so subdividing the edge does not change the component relation above rank $\alpha$. If exactly one endpoint has rank below $\alpha$, then every internal vertex of the replacement path has the rank of the endpoint of larger rank. After deleting the endpoint of rank below $\alpha$, the remaining subpath is attached to the component containing the endpoint of larger rank and to no other component. Finally, no edge of $I$ joins $\upset{t}\times\omega$ to $\upset{t'}\times\omega$ when $t,t'\in T^\alpha$ are distinct. The component description therefore follows from \cref{lem:components-inflation}.
\end{proof}

We shall use the following elementary degree observation.

\begin{remark}\label{rem:degree}
Let $T$ be pruned and let $I=\rayinfl{G}$.
\begin{enumerate}
  \item If $t$ is not the root and $n\ge1$, then $(t,n)$ has degree at least three: its two horizontal neighbours and a lower neighbour are distinct.
  \item If $u$ is a limit node and $m\in\omega$, then $(u,m)$ has degree at least three. Since $T$ is pruned, choose $w>_Tu$ and put $u^+=w\restr(\treeht_T(u)+1)$. Then $u^+$ is an immediate successor of $u$, and $(C_m(T,u),m)$, $(u,m+1)$, and $(u^+,m)$ are distinct neighbours of $(u,m)$.
\end{enumerate}
\end{remark}

\subsection{Ends and ray graphs}

The set of ends of a graph $G$ is denoted by $\Omega(G)$. Let $\mathcal R$ be a family of pairwise disjoint rays. An \emph{$\mathcal R$-path} is a finite path whose endvertices lie on distinct rays of $\mathcal R$ and whose internal vertices avoid all rays of $\mathcal R$. A family of $\mathcal R$-paths is \emph{independent} if any two of its members may intersect only in their endvertices. Given such a family $\mathcal P$, the graph associated with $(\mathcal R,\mathcal P)$ has vertex set $\mathcal R$, and two rays are adjacent if infinitely many paths in $\mathcal P$ join them.

A \emph{ray graph} for an end $\varepsilon$ is a connected graph associated with a family $\mathcal R\subseteq\varepsilon$ of $\degop(\varepsilon)$ pairwise disjoint rays and an independent family of $\mathcal R$-paths. For a cardinal $\kappa$, let $\HC(\kappa)$ denote the assertion that every end of degree $\kappa$ in every graph has a ray graph.

For an infinite cardinal $\kappa$, a \emph{$\kappa$-star of rays} consists of a centre ray $R$, pairwise disjoint leaf rays $R_i$ for $i<\kappa$, and, for each $i$, infinitely many pairwise disjoint $R_i$--$R$ paths, such that paths belonging to distinct leaves meet only on $R$.

\section{Club arguments}\label{sec:club}

A \emph{rank function} on a set $X$ is a map $\rho_X:X\to\omega_1$. Put
\[
  X^{<\alpha}:=\{x\in X:\rho_X(x)<\alpha\}.
\]
We say that $\rho_X$ has \emph{countable initial segments} if $X^{<\alpha}$ is countable for every $\alpha<\omega_1$.

\begin{lemma}\label{lem:club-sync}
Let $X$ and $Y$ carry rank functions with countable initial segments, and let $f:X\to Y$ be injective. There is a club $C\subseteq\omega_1$, consisting of non-zero limit ordinals, such that
\begin{equation}\label{eq:club-injection}
  f[X^{<\alpha}]=f[X]\cap Y^{<\alpha}
  \qquad\text{for every }\alpha\in C.
\end{equation}
If $f$ is bijective, then $f[X^{<\alpha}]=Y^{<\alpha}$ for every $\alpha\in C$.
\end{lemma}

\begin{proof}
For $\gamma<\omega_1$, define
\[
 a(\gamma):=\sup\{\rho_Y(f(x))+1:x\in X^{<\gamma}\}
\]
and
\[
 b(\gamma):=\sup\{\rho_X(x)+1:x\in X,\ \rho_Y(f(x))<\gamma\}.
\]
Both ordinals are below $\omega_1$. The first set is countable because $X^{<\gamma}$ is countable. For the second, injectivity of $f$ maps the set $\{x\in X:\rho_Y(f(x))<\gamma\}$ bijectively onto a subset of the countable set $Y^{<\gamma}$.

Let $C$ be a club of non-zero limit ordinals closed under both $a$ and $b$, in the sense that
\[
  \gamma<\alpha\in C
  \quad\Longrightarrow\quad
  a(\gamma)<\alpha\text{ and }b(\gamma)<\alpha.
\]
Fix $\alpha\in C$. If $x\in X^{<\alpha}$, choose $\gamma<\alpha$ with $\rho_X(x)<\gamma$. Then $\rho_Y(f(x))<a(\gamma)<\alpha$. Conversely, if $f(x)\in Y^{<\alpha}$, choose $\gamma<\alpha$ with $\rho_Y(f(x))<\gamma$. Then $\rho_X(x)<b(\gamma)<\alpha$. This proves \eqref{eq:club-injection}. If $f$ is bijective, then $f[X]=Y$, and the right-hand side of \eqref{eq:club-injection} is $Y^{<\alpha}$.
\end{proof}

Every rank-initial segment of a ray inflation over an $\omega_1$-tree is countable. If $J$ is a subdivision of such an inflation, then $J^{<\alpha}$ is also countable: below rank $\alpha$ there are only countably many original vertices and countably many edges, and each edge is replaced by a finite path. Hence \cref{lem:club-sync} applies both to isomorphisms involving subdivisions and to graph embeddings between ray inflations.

We first record the standard reflection fact for specialness.

\begin{lemma}[Abraham; see \protect\cite{DevlinJohnsbraten}]\label{lem:special-on-club}
Let $T$ be an $\omega_1$-tree and let $C\subseteq\omega_1$ be a club. If
\[
  T\restr C:=\bigcup_{\alpha\in C}T^\alpha
\]
is special with the inherited order, then $T$ is special.
\end{lemma}

For a sparse $T$-graph $G$ and $v\in T$, define
\begin{equation}\label{eq:Bv}
  B_G(v):=\bigcup_{u>_Tv}\bigl(N_G^-(u)\cap\sdownset{v}\bigr).
\end{equation}
Thus $B_G(v)$ consists of lower endpoints strictly below $v$ of edges whose upper endpoint lies strictly above $v$.

\begin{lemma}\label{lem:club-nonuniformity}
Let $T$ be a non-special $\omega_1$-tree and let $G$ be a sparse $T$-graph. For every club $C\subseteq\omega_1$, there are $\alpha\in C$ and $v\in T^\alpha$ such that $B_G(v)$ is infinite.
\end{lemma}

\begin{proof}
Suppose otherwise. The set $\{0\}\cup\acc(C)$ is club. Replacing $C$ by this set, we may assume that $0\in C$, that every non-zero member of $C$ is a limit point of $C$, and that $B_G(v)$ is finite for every $v\in T\restr C$. Enumerate
\[
  C=\{\gamma_\xi:\xi<\omega_1\}
\]
in increasing continuous order, with $\gamma_0=0$, and let $U=T\restr C$ with the inherited order.

We construct a sparse $U$-graph $H$. Join every successor node of $U$ to its predecessor in $U$. Let $u\in T^{\gamma_\lambda}$ be a limit node of $U$. For $x<_Tu$, define
\[
  \pi_u(x):=u\restr\gamma_{\eta(x)},
  \qquad
  \eta(x):=\min\{\eta<\lambda:\treeht_T(x)\le\gamma_\eta\}.
\]
The sequence $\langle\pi_u(C_n(T,u)):n\in\omega\rangle$ is nondecreasing in the tree order. Delete repeated terms. The remaining sequence is strictly increasing and cofinal in $\sdownset{u}\cap U$; declare its terms to be the lower neighbours of $u$ in $H$.

We claim that $H$ has uniformly finite adhesion. Fix a limit node $v\in U$. If $u>_Uv$ is a successor node of $U$, its predecessor in $U$ does not lie strictly below $v$. If $u$ is a limit node and a projected lower neighbour $\pi_u(C_n(T,u))$ lies below $v$, then $C_n(T,u)<_Tv$, so $C_n(T,u)\in B_G(v)$. Moreover, if
\[
  \pi_u(C_n(T,u))=u\restr\gamma_\eta<_Tv,
\]
then this node equals $v\restr\gamma_\eta$. For $b\in B_G(v)$, put
\[
  \eta(b):=\min\{\xi:\treeht_T(b)\le\gamma_\xi\}.
\]
Consequently, every lower neighbour below $v$ of every node above $v$ belongs to the finite set
\[
  S_v:=\{v\restr\gamma_{\eta(b)}:b\in B_G(v),\ \gamma_{\eta(b)}<\treeht_T(v)\}.
\]
Thus $H$ has uniformly finite adhesion. By \cref{thm:KP}, the tree $U$ is special, and \cref{lem:special-on-club} implies that $T$ is special, a contradiction.
\end{proof}

\section{Isomorphisms onto subdivisions}\label{sec:whole-graph}

We first derive the finite boundary property of uniformly finite adhesion that will be used in the proof. For $m\in\omega$, write $[m]:=\{0,\ldots,m-1\}$.

\begin{lemma}[Geschke, Kurkofka, Melcher, and Pitz~\protect-\cite{GeschkeEtAl}]\label{lem:uniform-boundary}
Let $G$ be a sparse $S$-graph with uniformly finite adhesion, let $s$ be a limit node of $S$, and let $F_s\subseteq\sdownset{s}$ be a finite witness at $s$. Suppose that two original vertices $(x,k)$ and $(y,\ell)$ are adjacent in $\rayinfl{G}$, with $x>_S s$ and $y<_S s$. Then
\[
  (y,\ell)\in F_s\times[|F_s|].
\]
\end{lemma}

\begin{theorem}\label{thm:intro-nonsubdivision}
Let $T$ and $S$ be pruned $\omega_1$-trees. Let $G_T$ be a sparse $T$-graph and let $G_S$ be a sparse $S$-graph with uniformly finite adhesion. If $T$ is not special, then $\rayinfl{G_T}$ is not isomorphic to a subdivision of $\rayinfl{G_S}$.
\end{theorem}
\begin{proof}
Put $I_T=\rayinfl{G_T}$ and $I_S=\rayinfl{G_S}$. Suppose that $J$ is a subdivision of $I_S$ and that
\[
  f:I_T\longrightarrow J
\]
is an isomorphism.

Equip $I_T$ and $J$ with the rank functions from \cref{sec:preliminaries}. By \cref{lem:club-sync}, there is a club $C\subseteq\omega_1$ of non-zero limit ordinals such that
\begin{equation}\label{eq:whole-sync}
  f[(I_T)^{<\alpha}]=J^{<\alpha}
  \qquad\text{for every }\alpha\in C.
\end{equation}
By \cref{lem:club-nonuniformity}, choose $\alpha\in C$ and $v\in T^\alpha$ such that $B_{G_T}(v)$ is infinite.

Equation~\eqref{eq:whole-sync} implies that $f$ restricts to an isomorphism
\[
  I_T-(I_T)^{<\alpha}\longrightarrow J-J^{<\alpha}.
\]
Therefore $f[K_{I_T}(v,\alpha)]$ is a component of $J-J^{<\alpha}$. By \cref{lem:components-subdivision}, there is a unique $s\in S^\alpha$ such that
\begin{equation}\label{eq:component-correspondence}
  f[K_{I_T}(v,\alpha)]=K_J(s,\alpha).
\end{equation}
The ordinal $\alpha$ is a non-zero limit, so $s$ is a limit node. Fix a finite witness $F_s\subseteq\sdownset{s}$ to uniformly finite adhesion at $s$.

\smallskip
\noindent\emph{Step 1: the image of a tail of $R_v$ is a tail of $R_s$.}
For $n\ge1$, put
\[
  p_n=(v,n),
  \qquad
  q_n=(C_n(T,v),n).
\]
The node $C_n(T,v)$ is non-root for $n\ge1$. Hence both $p_n$ and $q_n$ have degree at least three by \cref{rem:degree}, and their images are original vertices of $J$. Write
\[
  f(p_n)=(x_n,k_n),
  \qquad
  f(q_n)=(y_n,\ell_n).
\]
By \eqref{eq:component-correspondence}, $x_n\ge_S s$. Since $q_n\in(I_T)^{<\alpha}$, equation~\eqref{eq:whole-sync} gives $\treeht_S(y_n)<\alpha=\treeht_S(s)$.

The vertices $f(p_n)$ and $f(q_n)$ are adjacent original vertices of $J$. Hence they are the endpoints of an edge of $I_S$ whose replacement path in $J$ has length one. Suppose that $x_n>_S s$. Since the edge is not horizontal, its first coordinates $x_n$ and $y_n$ are comparable in $S$. The nodes $s$ and $y_n$ both belong to the well-ordered set $\downset{x_n}$. Since $\treeht_S(y_n)<\treeht_S(s)$, we have $y_n<_S s$. Thus \cref{lem:uniform-boundary} yields
\[
  f(q_n)\in F_s\times[|F_s|].
\]
The vertices $q_n$ are distinct and $f$ is injective, so this occurs for only finitely many $n$. Consequently, there is $N\ge1$ such that
\begin{equation}\label{eq:base-s}
  f(v,n)=(s,k_n)
  \qquad\text{for every }n\ge N.
\end{equation}

For $n\ge N$, the adjacent original vertices $f(v,n)$ and $f(v,n+1)$ have the same first coordinate. Hence they are joined in $I_S$ by an unsubdivided horizontal edge, and
\[
  |k_{n+1}-k_n|=1.
\]
The sequence $(k_n)_{n\ge N}$ is injective. The sign of $k_{n+1}-k_n$ cannot change: a change from $+1$ to $-1$, or from $-1$ to $+1$, would repeat the preceding value. Hence the sequence is either strictly increasing or strictly decreasing. It cannot be strictly decreasing for infinitely many terms in $\omega$, so it is strictly increasing. Consequently,
\begin{equation}\label{eq:tail-to-tail}
  f[\{(v,n):n\ge N\}]=\{(s,k):k\ge k_N\}.
\end{equation}

\smallskip
\noindent\emph{Step 2: the infinite set $B_{G_T}(v)$ contradicts uniformly finite adhesion at $s$.}
Recall that $v\in T^\alpha$, that $B_{G_T}(v)$ is infinite, and that $s\in S^\alpha$ is the unique node satisfying
\[
  f[K_{I_T}(v,\alpha)]=K_J(s,\alpha).
\]
Let $F_s\subseteq\sdownset{s}$ be the finite set witnessing uniformly finite adhesion of $G_S$ at $s$.

For $u>_Tv$, define
\[
  L_v(u):=N^-_{G_T}(u)\cap\sdownset{v}.
\]
Each set $L_v(u)$ is finite. If $u$ is a successor node, then its only lower neighbour is $\treepred_T(u)\ge_Tv$, so $L_v(u)=\varnothing$. If $u$ is a limit node, let
\[
  n_v(u):=\min\{n<\omega:v\le_TC_n(T,u)\}.
\]
The lower-neighbour sequence of $u$ is strictly increasing, and both $v$ and every $C_n(T,u)$ lie in the linearly ordered set $\sdownset{u}$. Therefore
\[
  L_v(u)=\{C_n(T,u):n<n_v(u)\},
\]
which is finite. Since
\[
  B_{G_T}(v)=\bigcup_{u>_Tv}L_v(u)
\]
is infinite, infinitely many distinct nodes above $v$ contribute lower neighbours below $v$.

Let $r_T$ be the root of $T$. We recursively choose pairs $(u_i,b_i)$, $i<\omega$, such that
\[
  u_i>_Tv,\qquad b_i<_Tv,\qquad b_i\in N^-_{G_T}(u_i),
\]
the nodes $u_i$ are pairwise distinct, and the nodes $b_i$ are pairwise distinct and different from $r_T$. Indeed, after choosing the pairs with index below $i$, the set
\[
  \{r_T\}\cup\{b_j:j<i\}\cup\bigcup_{j<i}L_v(u_j)
\]
is finite. Choose $b_i$ in $B_{G_T}(v)$ outside this set and then choose $u_i>_Tv$ with $b_i\in N^-_{G_T}(u_i)$. The choice of $b_i$ ensures that $u_i$ is new. Moreover, $u_i$ cannot be a successor node, because the immediate predecessor of a successor $u_i>_Tv$ is at least $v$. Thus $u_i$ is a limit node, and there is $m_i<\omega$ such that
\[
  b_i=C_{m_i}(T,u_i).
\]
Put
\[
  a_i:=(u_i,m_i),\qquad d_i:=(b_i,m_i).
\]
Then $a_id_i$ is an edge of $I_T$. Both endpoints have degree at least three. For $a_i$, an immediate successor $u_i^+$ of $u_i$ gives the three distinct neighbours
\[
  d_i,\qquad (u_i,m_i+1),\qquad (u_i^+,m_i).
\]
For $d_i$, let $c_i=\treepred_T(b_i)$ when $b_i$ is a successor node and let $c_i=C_{m_i}(T,b_i)$ when $b_i$ is a limit node. Since $b_i\ne r_T$, the vertices
\[
  a_i,\qquad (b_i,m_i+1),\qquad (c_i,m_i)
\]
are three distinct neighbours of $d_i$. Hence $f(a_i)$ and $f(d_i)$ are original vertices of the subdivision $J$.

Write
\[
  f(a_i)=(x_i,r_i),\qquad f(d_i)=(y_i,\ell_i).
\]
Since $a_i\in K_{I_T}(v,\alpha)$, equation~\eqref{eq:component-correspondence} gives $x_i\ge_Ss$. Since $d_i\in(I_T)^{<\alpha}$, equation~\eqref{eq:whole-sync} gives $\treeht_S(y_i)<\alpha=\treeht_S(s)$. The adjacent vertices $f(a_i)$ and $f(d_i)$ are original vertices of $J$, so they are the endpoints of an unsubdivided edge of $I_S$. This edge is not horizontal, and therefore $x_i$ and $y_i$ are comparable. As $s\le_Sx_i$ and both $s$ and $y_i$ are predecessors of $x_i$, the inequality of their heights implies
\[
  y_i<_Ss.
\]

Only finitely many indices satisfy $x_i=s$. Indeed, if $x_i=s$, then $f(a_i)$ lies on $R_s$. By \eqref{eq:tail-to-tail}, every vertex $(s,k)$ with $k\ge k_N$ is the image of a vertex of the tail $\{(v,n):n\ge N\}$ of $R_v$. Since $u_i>_Tv$, the vertex $a_i$ does not lie on $R_v$, and injectivity of $f$ therefore forces
\[
  f(a_i)\in\{(s,k):k<k_N\}.
\]
This set is finite, whereas the vertices $f(a_i)$ are pairwise distinct. After passing to an infinite set of indices, we may consequently assume that $x_i>_Ss$ for every $i$. The edge of $I_S$ with endpoints $f(a_i)$ and $f(d_i)$ then has one endpoint above $s$ and the other below $s$. By \cref{lem:uniform-boundary},
\[
  f(d_i)\in F_s\times[|F_s|]
\]
for every remaining $i$. This is impossible, because the vertices $d_i$ are pairwise distinct, $f$ is injective, and $F_s\times[|F_s|]$ is finite.
\end{proof}

\section{Graph embeddings}\label{sec:subgraph}

Let $S$ and $T$ be $\omega_1$-trees, let $G_S$ and $G_T$ be sparse graphs on them, and write
\[
  I_S=\rayinfl{G_S},
  \qquad
  I_T=\rayinfl{G_T}.
\]
Fix a graph embedding
\[
  f:I_S\longrightarrow I_T;
\]
that is, $f$ is injective and preserves adjacency. It need not preserve non-adjacency.

By \cref{lem:club-sync}, there is a club $C\subseteq\omega_1$ of non-zero limit ordinals such that
\begin{equation}\label{eq:embedding-sync}
  f[(I_S)^{<\alpha}]=f[V(I_S)]\cap(I_T)^{<\alpha}
  \qquad\text{for every }\alpha\in C.
\end{equation}

\begin{lemma}\label{lem:cone-map}
For every $\alpha\in C$ and $s\in S^\alpha$, there is a unique node $\pi(s)\in T^\alpha$ such that
\begin{equation}\label{eq:cone-map}
  f[K_{I_S}(s,\alpha)]\subseteq K_{I_T}(\pi(s),\alpha).
\end{equation}
Moreover, if $\alpha<\beta$ belong to $C$, $s\in S^\alpha$, $s'\in S^\beta$, and $s<_S s'$, then
\begin{equation}\label{eq:coherence}
  \pi(s)=\pi(s')\restr\alpha.
\end{equation}
\end{lemma}

\begin{proof}
Every vertex of $K_{I_S}(s,\alpha)$ has rank at least $\alpha$ in $I_S$. Equation~\eqref{eq:embedding-sync} therefore implies that every vertex of its image has rank at least $\alpha$ in $I_T$. Since $K_{I_S}(s,\alpha)$ is connected, its image is contained in a unique component of $I_T-(I_T)^{<\alpha}$. By \cref{lem:components-inflation}, this component is $K_{I_T}(t,\alpha)$ for a unique $t\in T^\alpha$. Define $\pi(s)=t$.

Suppose that $s<_S s'$. Then
\[
  K_{I_S}(s',\beta)\subseteq K_{I_S}(s,\alpha),
\]
so
\[
  f[K_{I_S}(s',\beta)]
  \subseteq K_{I_T}(\pi(s'),\beta)\cap K_{I_T}(\pi(s),\alpha).
\]
The set on the left is nonempty. Moreover,
\[
  K_{I_T}(\pi(s'),\beta)
  \subseteq K_{I_T}(\pi(s')\restr\alpha,\alpha).
\]
The components of $I_T-(I_T)^{<\alpha}$ are pairwise disjoint. Hence
\[
  K_{I_T}(\pi(s),\alpha)
  =K_{I_T}(\pi(s')\restr\alpha,\alpha),
\]
and the uniqueness in \cref{lem:components-inflation} gives
\[
  \pi(s)=\pi(s')\restr\alpha.
\]
In particular, $\pi(s)<_T\pi(s')$.
\end{proof}

Fix a non-zero limit ordinal $\alpha<\omega_1$ and $t\in T^\alpha$. Define the set of vertices of the component $K_{I_T}(t,\alpha)$ that have a neighbour of rank below $\alpha$ by
\[
  X_\alpha(t):=
  \{x\in V(K_{I_T}(t,\alpha)):
  N_{I_T}(x)\cap(I_T)^{<\alpha}\ne\varnothing\}.
\]

\begin{lemma}\label{lem:lower-boundary}
The connected component of $I_T[X_\alpha(t)]$ containing $R_t$ is exactly
$R_t$. Every other non-empty connected component of $I_T[X_\alpha(t)]$ is
a finite initial segment of $R_u$ for some limit node $u>_Tt$. Consequently,
$R_t$ is the unique connected component of $I_T[X_\alpha(t)]$ that contains
a ray.
\end{lemma}

\begin{proof}
We first determine the intersection of $X_\alpha(t)$ with each horizontal
ray contained in $K_{I_T}(t,\alpha)$.

Since $\alpha$ is a limit ordinal, the node $t\in T^\alpha$ is a limit node.
For every $n\in\omega$, the lower neighbour $C_n(T,t)$ has height below
$\alpha$, and the edge
\[
  (t,n)(C_n(T,t),n)
\]
shows that $(t,n)\in X_\alpha(t)$. Therefore
\[
  X_\alpha(t)\cap R_t=V(R_t).
\]

Now fix $u>_Tt$. Suppose first that $u$ is a successor node. Its only lower
neighbour in $G_T$ is $\treepred_T(u)$. Since $t<_Tu$, the immediate
predecessor of $u$ satisfies $t\le_T\treepred_T(u)$, and hence has height at
least $\alpha$. Thus no vertex of $R_u$ has a neighbour of rank below
$\alpha$, so
\[
  X_\alpha(t)\cap R_u=\varnothing.
\]

Suppose next that $u$ is a limit node. Its lower neighbours form the strictly
increasing cofinal sequence
\[
  C_0(T,u)<_TC_1(T,u)<_T\cdots<_Tu.
\]
Since $t<_Tu$, there is a least integer
\[
  k(u,t):=\min\{n<\omega:t\le_TC_n(T,u)\}.
\]
For $n<k(u,t)$, the nodes $C_n(T,u)$ and $t$ both belong to the linearly
ordered set $\sdownset{u}$, and the minimality of $k(u,t)$ gives
$C_n(T,u)<_Tt$. Hence $(u,n)$ is adjacent to a vertex of rank below
$\alpha$. For $n\ge k(u,t)$, we have
\[
  t\le_TC_{k(u,t)}(T,u)\le_TC_n(T,u),
\]
so the unique lower neighbour of $(u,n)$ has rank at least $\alpha$.
All horizontal neighbours of $(u,n)$ have rank $\treeht_T(u)>\alpha$, and
all neighbours arising from nodes above $u$ also have rank greater than
$\alpha$. Therefore $(u,n)$ has no neighbour of rank below $\alpha$.
Consequently,
\[
  X_\alpha(t)\cap R_u
  =\{(u,n):n<k(u,t)\},
\]
which is a finite initial segment of $R_u$.

It remains to verify that no edge between two distinct horizontal rays has
both endpoints in $X_\alpha(t)$. Consider such an edge
\[
  (v,n)(w,n),
  \qquad v<_Tw.
\]
If both endpoints belong to $K_{I_T}(t,\alpha)$, then $t\le_Tv$. At
coordinate $n$, the vertex $(v,n)$ is the unique lower neighbour of
$(w,n)$. Since $\treeht_T(v)\ge\alpha$, this lower neighbour does not have
rank below $\alpha$; the preceding description of the neighbours of
$(w,n)$ shows that none of its other neighbours has rank below $\alpha$
either. Hence $(w,n)\notin X_\alpha(t)$. Thus no non-horizontal edge has
both endpoints in $X_\alpha(t)$.

The induced graph $I_T[X_\alpha(t)]$ is therefore the disjoint union of the
whole horizontal ray $R_t$ and the finite initial segments described above.
The asserted component structure follows.
\end{proof}

\begin{proposition}\label{prop:horizontal-rigidity}
For every $\alpha\in C$ and $s\in S^\alpha$, the image $f[R_s]$ is a tail of
$R_{\pi(s)}$.
\end{proposition}

\begin{proof}
Fix $\alpha\in C$ and $s\in S^\alpha$, and put $t=\pi(s)$. Since $\alpha$ is a non-zero limit ordinal, $s$ is a limit node. For each $n<\omega$, set
\[
  x_n:=f(s,n),\qquad y_n:=f(C_n(S,s),n).
\]
The vertices $(s,n)$ and $(C_n(S,s),n)$ are adjacent in $I_S$. Moreover, $(C_n(S,s),n)$ has rank below $\alpha$, whereas $(s,n)$ belongs to $K_{I_S}(s,\alpha)$. Equations~\eqref{eq:embedding-sync} and~\eqref{eq:cone-map} therefore give
\[
  y_n\in(I_T)^{<\alpha},\qquad x_n\in K_{I_T}(t,\alpha).
\]
Because $f$ preserves adjacency, $x_n$ has a neighbour of rank below $\alpha$, and hence $x_n\in X_\alpha(t)$.

The vertices $x_n$ are pairwise distinct and consecutive vertices are adjacent, so they form a ray in $I_T[X_\alpha(t)]$. By \cref{lem:lower-boundary}, this ray lies in $R_t$. Write $x_n=(t,k_n)$. Adjacency on $R_t$ gives $|k_{n+1}-k_n|=1$, and injectivity gives $k_n\ne k_m$ whenever $n\ne m$. The signs of two consecutive differences cannot be opposite, since that would imply $k_{n+2}=k_n$; thus all differences have the same sign. They cannot all equal $-1$, because $\omega$ has no infinite strictly decreasing sequence. Therefore $k_{n+1}=k_n+1$ for every $n$, so $k_n=k_0+n$ and
\[
  f[R_s]=\{(t,k):k\ge k_0\},
\]
a tail of $R_t=R_{\pi(s)}$.
\end{proof}

\begin{lemma}\label{lem:level-injective}
For every $\alpha\in C$, the map $\pi\restr S^\alpha:S^\alpha\to T^\alpha$ is injective.
\end{lemma}

\begin{proof}
Suppose that $s,s'\in S^\alpha$ are distinct and $\pi(s)=\pi(s')=t$. The rays $R_s$ and $R_{s'}$ are disjoint, and injectivity of $f$ implies that their images are disjoint. By \cref{prop:horizontal-rigidity}, both images are tails of $R_t$, but two tails of the same ray meet. This is a contradiction.
\end{proof}

\begin{theorem}\label{thm:intro-rigidity}
Let $S$ and $T$ be $\omega_1$-trees, and let $G_S$ and $G_T$ be sparse graphs on $S$ and $T$, respectively. If
\[
  f:\rayinfl{G_S}\longrightarrow\rayinfl{G_T}
\]
is a graph embedding, then there are a club $C\subseteq\omega_1$, consisting of non-zero limit ordinals, and a level-preserving order embedding
\[
  \pi:S\restr C\longrightarrow T\restr C
\]
such that, for every $\alpha\in C$ and $s\in S^\alpha$, the image $f[R_s]$ is a tail of $R_{\pi(s)}$.
\end{theorem}

\begin{proof}
The map $\pi$ preserves levels and strict order by \cref{lem:cone-map}. It remains to prove that it reflects order. Suppose that $s\in S^\alpha$, $s'\in S^\beta$, $\alpha<\beta$ belong to $C$, and $\pi(s)<_T\pi(s')$. Let $r=s'\restr\alpha$. By \eqref{eq:coherence},
\[
  \pi(r)=\pi(s')\restr\alpha.
\]
The node $\pi(s)$ also lies at level $\alpha$ below $\pi(s')$, so uniqueness of the predecessor at level $\alpha$ gives $\pi(s)=\pi(r)$. By \cref{lem:level-injective}, $s=r$, and hence $s<_S s'$. The statement about horizontal rays is \cref{prop:horizontal-rigidity}.
\end{proof}

\begin{corollary}\label{cor:intro-almost-suslin}
Let $T$ be an almost-Suslin $\omega_1$-tree and let $S$ be a special $\omega_1$-tree. For arbitrary sparse graphs $G_T$ on $T$ and $G_S$ on $S$, the graph $\rayinfl{G_T}$ contains no subgraph isomorphic to $\rayinfl{G_S}$.
\end{corollary}

\begin{proof}
Suppose that such a subgraph exists, and let $f:\rayinfl{G_S}\to\rayinfl{G_T}$ be the corresponding embedding. By \cref{thm:intro-rigidity}, there are a club $C\subseteq\omega_1$ and a level-preserving order embedding
\[
  \pi:S\restr C\longrightarrow T\restr C.
\]
Write
\[
  S=\bigcup_{n<\omega}A_n,
\]
where every $A_n$ is an antichain. Since $\pi$ reflects comparability, each set
\[
  B_n:=\pi[A_n\cap(S\restr C)]
\]
is an antichain in $T$.

For every $\alpha\in C$, the level $S^\alpha$ is non-empty and $\pi$ preserves levels. Hence
\[
  C\subseteq\bigcup_{n<\omega}\{\treeht_T(t):t\in B_n\}.
\]
For $n<\omega$, let
\[
  H_n:=\{\treeht_T(t):t\in B_n\}.
\]
If every $H_n$ were nonstationary, then their countable union would be nonstationary, because the nonstationary ideal on $\omega_1$ is countably complete. This is impossible, since the union contains the club $C$. Hence $H_n$ is stationary for some $n$, and $B_n$ is an antichain in $T$ whose set of heights is stationary. This contradicts the almost-Suslin property of $T$.
\end{proof}

\begin{remark}\label{rem:subgraph-strength}
The embedding in \cref{cor:intro-almost-suslin} need only preserve adjacency; it is not required to preserve non-adjacency. The proof of \cref{thm:intro-rigidity} does not require either tree to be pruned, Aronszajn, semi-special, or non-special.
\end{remark}

\section{Ends of finite-adhesion inflations}\label{sec:ends}

We first show that finite adhesion excludes uncountable stars of rays.

\begin{proposition}\label{prop:no-star}
Let $T$ be a pruned Aronszajn tree and let $G$ be a sparse $T$-graph with finite adhesion. Then $\rayinfl{G}$ contains no $\aleph_1$-star of rays.
\end{proposition}

\begin{proof}
The argument is an adaptation of the proof of Geschke, Kurkofka, Melcher, and Pitz~\cite[Theorem~6.6]{GeschkeEtAl}. In their proof, the finite boundary property for the ray inflation is supplied by~\cite[Lemma~6.5]{GeschkeEtAl}, which is derived from property $(\star)$ in~\cite[Proposition~6.4]{GeschkeEtAl}. Here the finite-adhesion hypothesis provides the corresponding finite boundary for successor cones, and the same counting argument can be applied.

Put $I=\rayinfl{G}$ and suppose that $I$ contains an $\aleph_1$-star of rays. Let $\mathcal S$ denote the subgraph formed by its centre ray $R$, its leaf rays $R_i$ for $i<\omega_1$, and all the specified connecting paths. For each $i<\omega_1$, let $D_i$ be the component of $\mathcal S-R$ containing $R_i$. The subgraphs $D_i$ are pairwise vertex-disjoint, and for each $i$ there are infinitely many pairwise disjoint $D_i$--$R$ paths in $\mathcal S$.

The ray $R$ is countable, so there is a successor ordinal $\alpha<\omega_1$ such that
\[
  V(R)\subseteq T^{<\alpha}\times\omega.
\]
Since $T^{<\alpha}\times\omega$ is countable and the subgraphs $D_i$ are pairwise vertex-disjoint, the set
\[
  J_0:=\{i<\omega_1:V(D_i)\cap(T^{<\alpha}\times\omega)\ne\varnothing\}
\]
is countable. For every $i\in\omega_1\setminus J_0$, the connected graph $D_i$ is contained in one component of $I-I^{<\alpha}$. Since $T^\alpha$ is countable, \cref{lem:components-inflation} yields an uncountable set $I_0\subseteq\omega_1\setminus J_0$ and a node $t\in T^\alpha$ such that
\[
  D_i\subseteq I[\upset{t}\times\omega]
  \qquad(i\in I_0).
\]
Because $R_t$ is countable and the $D_i$, $i\in I_0$, are pairwise vertex-disjoint, the set
\[
  J_1:=\{i\in I_0:V(D_i)\cap V(R_t)\ne\varnothing\}
\]
is countable. For $i\in I_0\setminus J_1$, every vertex of $D_i$ has first coordinate strictly above $t$ and therefore rank at least $\alpha+1$. Applying \cref{lem:components-inflation} at level $\alpha+1$, choose a node $p\in T^{\alpha+1}$ above $t$ and an uncountable set $I_1\subseteq I_0\setminus J_1$ such that
\[
  D_i\subseteq I[\upset{p}\times\omega]
  \qquad(i\in I_1).
\]

Because $p\in T^{\alpha+1}$ lies above $t\in T^\alpha$, the node $p$ is a successor node and
\[
  \treepred_T(p)=t.
\]
At this point we use the successor-cone characterisation of finite adhesion, not merely its defining formulation in terms of components above limit levels. By Kurkofka and Pitz~\cite[Lemma~4.2]{KurkofkaPitz}, a $T$-graph has finite adhesion if and only if
\[
  N_G(\upset{q})\text{ is finite for every successor node }q\in T.
\]
Applying this characterisation to $q=p$, the set
\[
  A:=N_G(\upset{p})
\]
is finite. Every edge with one endpoint in $\upset{p}\times\omega$ and the other endpoint outside $\upset{p}\times\omega$ has its outside endpoint in
\begin{equation}\label{eq:finite-adhesion-boundary}
  (\{t\}\times\omega)\cup(A\times[|A|]).
\end{equation}
Indeed, the edges joining $(p,n)$ to $(t,n)$, for $n\in\omega$, have their outside endpoints in $\{t\}\times\omega$. Every other such edge has the form
\[
  (u,n)(C_n(T,u),n),
\]
where $u\ge_Tp$ is a limit node and $C_n(T,u)\notin\upset{p}$. Then $C_n(T,u)\in A$. For a fixed $u$, the terms of the increasing sequence $\langle C_m(T,u):m\in\omega\rangle$ that lie outside $\upset{p}$ form an initial segment of distinct elements of the finite set $A$. Hence $n<|A|$, which proves \eqref{eq:finite-adhesion-boundary}.

Fix $i\in I_1$. Every $D_i$--$R$ path must leave $\upset{p}\times\omega$, because $D_i\subseteq I[\upset{p}\times\omega]$ while $R\subseteq T^{<\alpha}\times\omega$. By \eqref{eq:finite-adhesion-boundary}, such a path meets either $\{t\}\times\omega$ or $A\times[|A|]$. The paths belonging to the fixed leaf $i$ are pairwise disjoint, so only finitely many of them can meet the finite set $A\times[|A|]$. Consequently, infinitely many meet $R_t=\{t\}\times\omega$, and they meet it in infinitely many distinct vertices.

If $i,j\in I_1$ are distinct, paths belonging to leaf $i$ and paths belonging to leaf $j$ can meet only on the centre ray $R$. Since $R\cap R_t=\varnothing$, the subsets of $R_t$ used by the two path systems are disjoint. We have therefore obtained uncountably many pairwise disjoint nonempty subsets of the countable set $V(R_t)$, a contradiction.
\end{proof}

A connected ray graph on $\aleph_1$ rays yields an $\aleph_1$-star of rays~\cite[Lemma~3.1]{GeschkeEtAl}. We can now prove the following theorem.

\begin{theorem}\label{thm:intro-end}
Let $T$ be a pruned semi-special $\omega_1$-tree, and let $G$ be a sparse $T$-graph with finite adhesion. Then $\rayinfl{G}$ has exactly one end, this end has degree $\aleph_1$, and it has no ray graph.
\end{theorem}

\begin{proof}
The tree $T$ is Aronszajn and has cardinality $\aleph_1$. By \cite[Lemma~6.2]{GeschkeEtAl}, the ray inflation $\rayinfl{G}$ has exactly one end and its degree is $|T|=\aleph_1$. If this end had a ray graph, \cite[Lemma~3.1]{GeschkeEtAl} would yield an $\aleph_1$-star of rays, contrary to \cref{prop:no-star}.
\end{proof}

Under $\diamondsuit_{\omega_1}$ there exists a pruned semi-special non-special Aronszajn tree; see Baumgartner's announcement~\cite{Baumgartner} and Hanazawa's construction~\cite{Hanazawa}. This assertion does not include the almost-Suslin property. For that strengthening we use $\diamondsuit_{\omega_1}^{*}$ and \cref{cor:diamond-star-tree}.

\begin{theorem}\label{thm:diamond-example}
Assume $\diamondsuit_{\omega_1}$. There is a graph $I$ with exactly one end $\varepsilon$ such that
\begin{enumerate}
  \item $\degop(\varepsilon)=\aleph_1$ and $\varepsilon$ has no ray graph;
  \item if $S$ is a pruned $\omega_1$-tree and $G_S$ is a sparse $S$-graph with uniformly finite adhesion, then $I$ is not isomorphic to a subdivision of $\rayinfl{G_S}$.
\end{enumerate}
\end{theorem}

\begin{proof}
Choose a pruned semi-special non-special $\omega_1$-tree $T$ and a standard sparse $T$-graph $G_T$ with finite adhesion. Put $I=\rayinfl{G_T}$. Assertion (i) follows from \cref{thm:intro-end}, and (ii) follows from \cref{thm:intro-nonsubdivision}.
\end{proof}

In particular, the graph in \cref{thm:diamond-example} is not isomorphic to a subdivision of any standard uniformly finite-adhesion ray inflation associated with a special $\omega_1$-tree. The conclusion concerns isomorphism of the entire graph; it does not exclude a subdivision of such a ray inflation occurring as a subgraph.

The stronger combinatorial principle also permits the exclusion of copies without edge subdivision.

\begin{theorem}\label{thm:diamond-star-example}
Assume $\diamondsuit_{\omega_1}^{*}$. There is a graph $I$ with exactly one end $\varepsilon$ such that
\begin{enumerate}
  \item $\degop(\varepsilon)=\aleph_1$ and $\varepsilon$ has no ray graph;
  \item if $S$ is a special $\omega_1$-tree and $G_S$ is any sparse $S$-graph, then $I$ contains no subgraph isomorphic to $\rayinfl{G_S}$;
  \item if $S$ is a pruned $\omega_1$-tree and $G_S$ is a sparse $S$-graph with uniformly finite adhesion, then $I$ is not isomorphic to a subdivision of $\rayinfl{G_S}$.
\end{enumerate}
\end{theorem}

\begin{proof}
Choose a pruned semi-special non-special almost-Suslin Aronszajn tree $T$ by \cref{cor:diamond-star-tree}, let $G_T$ be a standard sparse $T$-graph with finite adhesion, and put $I=\rayinfl{G_T}$. Assertion (i) follows from \cref{thm:intro-end}. Assertion (ii) follows from \cref{cor:intro-almost-suslin}, and assertion (iii) follows from \cref{thm:intro-nonsubdivision}.
\end{proof}

\begin{remark}\label{rem:diamond-star-scope}
The conclusion in \cref{thm:diamond-star-example}(ii) excludes arbitrary subgraphs isomorphic to $\rayinfl{G_S}$; the subgraph need not be induced. It does not exclude a subdivision of $\rayinfl{G_S}$ occurring as a subgraph. This is the problem formulated in \cref{conj:almost-suslin-subdivision}.
\end{remark}

We now introduce a local property weaker than the existence of a ray graph.

\begin{definition}[The $\HCstar$ property]\label{def:HCstar}
Let $G$ be a graph and let $\varepsilon\in\Omega(G)$ be an end of degree $\kappa$. We say that $\varepsilon$ has property $\HCstar(\kappa)$ if there are a family $\mathcal R\subseteq\varepsilon$ of $\kappa$ pairwise disjoint rays and an independent family $\mathcal P$ of $\mathcal R$-paths such that the following holds: for every countable subfamily $\mathcal R'\subseteq\mathcal R$, there are countable subfamilies $\mathcal R''\subseteq\mathcal R$ and $\mathcal P''\subseteq\mathcal P$ such that
\[
  \mathcal R'\subseteq\mathcal R''
\]
and the graph associated with $(\mathcal R'',\mathcal P'')$ is connected.
\end{definition}

\begin{lemma}\label{lem:HC-implies-HCstar}
If $\HC(\kappa)$ holds, then every end of degree $\kappa$ has property $\HCstar(\kappa)$. In particular, if $\kappa$ is infinite, then every end of degree $\kappa$ containing a $\kappa$-star of rays has property $\HCstar(\kappa)$.
\end{lemma}

\begin{proof}
Assume first that $\HC(\kappa)$ holds, and let $\varepsilon$ be an end of degree $\kappa$ in a graph $G$. Choose a ray graph for $\varepsilon$, witnessed by a family $\mathcal R\subseteq\varepsilon$ of $\kappa$ pairwise disjoint rays and an independent family $\mathcal P$ of $\mathcal R$-paths. Let $H$ be the connected graph associated with $(\mathcal R,\mathcal P)$.

Fix a countable subfamily $\mathcal R'\subseteq\mathcal R$. There is a countable connected subgraph $H'\subseteq H$ containing $\mathcal R'$: enumerate $\mathcal R'$ and join each of its members to a fixed vertex of $H$ by a finite path. Put $\mathcal R''=V(H')$. For every edge $RR'$ of $H'$, choose countably infinitely many paths from $\mathcal P$ joining $R$ to $R'$, and let $\mathcal P''$ be their union. The families $\mathcal R''$ and $\mathcal P''$ are countable, $\mathcal P''$ is independent, and the graph associated with $(\mathcal R'',\mathcal P'')$ contains $H'$. Thus it is connected, and $(\mathcal R,\mathcal P)$ witnesses $\HCstar(\kappa)$.

For the final assertion, let $R$ be the centre ray of a $\kappa$-star and let $R_i$, $i<\kappa$, be its leaf rays. Trim every connecting path at its first intersection with $R$, and choose countably infinitely many paths for each leaf. Their union is an independent family $\mathcal P$ of paths for
\[
  \mathcal R:=\{R\}\cup\{R_i:i<\kappa\}.
\]
Given a countable $\mathcal R'\subseteq\mathcal R$, let $\mathcal R''$ consist of $R$ together with the leaf rays in $\mathcal R'$, and let $\mathcal P''$ consist of the chosen paths for those leaves. The associated graph is a connected star.
\end{proof}

\section{Open problems}\label{sec:open}

The preceding lemma motivates the following weakening of Halin's end degree conjecture.

\begin{conjecture}\label{conj:HCstar}
For every cardinal $\kappa$, every end of degree $\kappa$ has property $\HCstar(\kappa)$.
\end{conjecture}

We also record questions about the form of counterexamples to $\HC(\aleph_1)$.

\begin{question}\label{question:non-inflation}
\begin{enumerate}[label=\textup{(\alph*)}]
  \item Is there a graph which is not isomorphic to the ray inflation of any $T$-graph and witnesses the failure of $\HC(\aleph_1)$?
  \item Is it consistent with $\mathrm{ZFC}$ that such a graph exists?
  \item Is it consistent with $\mathrm{ZFC}$ that every counterexample to $\HC(\aleph_1)$ is isomorphic to the ray inflation of a $T$-graph on a special $\omega_1$-tree?
\end{enumerate}
\end{question}

Semi-specialness and non-specialness alone do not prevent a ray inflation over $T$ from containing a ray inflation over a special tree. Let $T_0$ be a pruned semi-special non-special $\omega_1$-tree and let $S$ be a pruned special $\omega_1$-tree. Form a new order tree $T$ by identifying the roots of $T_0$ and $S$ and declaring every node of $T_0$ other than the root incomparable with every node of $S$ other than the root. Then $T$ is pruned, semi-special, and non-special.

Choose partitions
\[
  \Succ(T_0)=\bigcup_{n<\omega}A_n^0,
  \qquad
  \Succ(S)=\bigcup_{n<\omega}A_n^1
\]
into antichains. For each $n<\omega$, the set $A_n^0\cup A_n^1$ is an antichain in $T$, and these sets partition $\Succ(T)$. Apply \cref{prop:standard-graph} to this partition and let $G_T$ be the resulting sparse $T$-graph with finite adhesion. If $u$ is a limit node of $S$, then every node below $u$ belongs to $S$, so the lower-neighbour sequence assigned to $u$ is contained in $S$. Therefore the restriction $G_S:=G_T[S]$ is a sparse $S$-graph, and the subgraph of $\rayinfl{G_T}$ induced by $S\times\omega$ is exactly $\rayinfl{G_S}$.

The almost-Suslin hypothesis excludes this construction at the level of ordinary subgraphs: \cref{cor:intro-almost-suslin} applies to arbitrary sparse graphs and excludes every subgraph isomorphic to a ray inflation over a special $\omega_1$-tree. By \cref{cor:diamond-star-tree}, $\diamondsuit_{\omega_1}^{*}$ implies the existence of pruned semi-special non-special almost-Suslin Aronszajn trees. The remaining problem asks whether the same conclusion holds for subdivisions.

\begin{conjecture}\label{conj:almost-suslin-subdivision}
Let $T$ be a pruned semi-special non-special almost-Suslin $\omega_1$-tree, and let $G_T$ be a standard sparse $T$-graph with finite adhesion. If $S$ is a pruned special $\omega_1$-tree and $G_S$ is a sparse $S$-graph, then $\rayinfl{G_T}$ contains no subgraph isomorphic to a subdivision of $\rayinfl{G_S}$.
\end{conjecture}

The proof of \cref{thm:intro-rigidity} does not extend directly to subdivisions. Let $s\in S^\alpha$ be a limit node. In an ordinary graph embedding, the edge
\[
  (s,n)(C_n(S,s),n)
\]
is mapped to an edge whose second endpoint has rank below $\alpha$. Hence the image of $(s,n)$ belongs to the set $X_\alpha(\pi(s))$, and \cref{lem:lower-boundary} controls the image of the entire horizontal ray $R_s$.

In a subdivision, the same source edge is represented by a finite path. The endpoints of this path correspond to $(s,n)$ and $(C_n(S,s),n)$, but the first vertex of the path whose rank is below $\alpha$ may be an internal subdivision vertex. Thus the original vertex corresponding to $(s,n)$ need not belong to $X_\alpha(\pi(s))$. The structure theorem for $I_T[X_\alpha(\pi(s))]$ then gives no immediate information about the subdivided image of $R_s$. A proof of \cref{conj:almost-suslin-subdivision} must therefore analyse the replacement paths that meet both rank regions, rather than only the original vertices of the subdivision.

\section*{Acknowledgements}

The first author was supported by Funda\c{c}\~ao de Amparo \`a Pesquisa do Estado de S\~ao Paulo (FAPESP), grant 2023/00595-6. The second author was supported by Funda\c{c}\~ao de Amparo \`a Pesquisa do Estado de S\~ao Paulo (FAPESP), grant 2025/09425-1. The third author was supported by Conselho Nacional de Desenvolvimento Cient\'ifico e Tecnol\'ogico (CNPq), grant 165761/2021-0.

\end{document}